# An Unsupervised Deep Learning Approach for Scenario Forecasts

Yize Chen*, Xiyu Wang†, and Baosen Zhang*
*Department of Electrical Engineering, University of Washington, Seattle, WA, USA
{yizechen, zhangbao}@uw.edu
†Department of Electrical Engineering, Tsinghua University, Beijing, China
{xiyu-wan14}@mails.tsinghua.edu.cn

***Abstract*—In this paper, we propose a novel scenario forecasts approach which can be applied to a broad range of power system operations (e.g., wind, solar, load) over various forecasts horizons and prediction intervals. This approach is model-free and data-driven, producing a set of scenarios that represent possible future behaviors based only on historical observations and point forecasts. It first applies a newly-developed unsupervised deep learning framework, the generative adversarial networks, to learn the intrinsic patterns in historical renewable generation data. Then by solving an optimization problem, we are able to quickly generate large number of realistic future scenarios. The proposed method has been applied to a wind power generation and forecasting dataset from national renewable energy laboratory. Simulation results indicate our method is able to generate scenarios that capture spatial and temporal correlations. Our code and simulation datasets are freely available online.**

## I. Introduction

The integration of high penetration of renewable generation into power systems calls for a growing need to model the uncertain and intermittent characteristics of these resources. An important method used in characterizing the behavior of renewable resources is *scenario generation*, where a set of possible future realizations are provided for the system operator. Compared to deterministic point forecasts or probabilistic forecasts [1], scenario forecasts could not only inform users of the uncertainty about the future, but also reflect the temporal dependence of renewable power generation [2], [3]. The information provided by these generated scenarios is valuable for a host of decision-making and stochastic optimization problems, such as the economic dispatch of renewables [4], [5], unit commitment [6], [7] and many others. Therefore, in recent years, many algorithms have been introduced for various applications, from load forecasting to wind and solar power generations.

One of the biggest challenges of scenario forecasts is the difficulty of modeling and learning the underlying stochastic processes that drives renewable power generation [8]. Previous statistical or physical methods like first-order autoregressive model [9], ensemble methods and Gaussian Copula [2], [10], [11] either required strong statistical assumptions or detailed physical measurements and modeling. What's more, most of these methods focus on capturing the marginal distribution of each individual time slots of the forecasting horizon, while paying less attention to the temporal correlations in the scenarios [12].

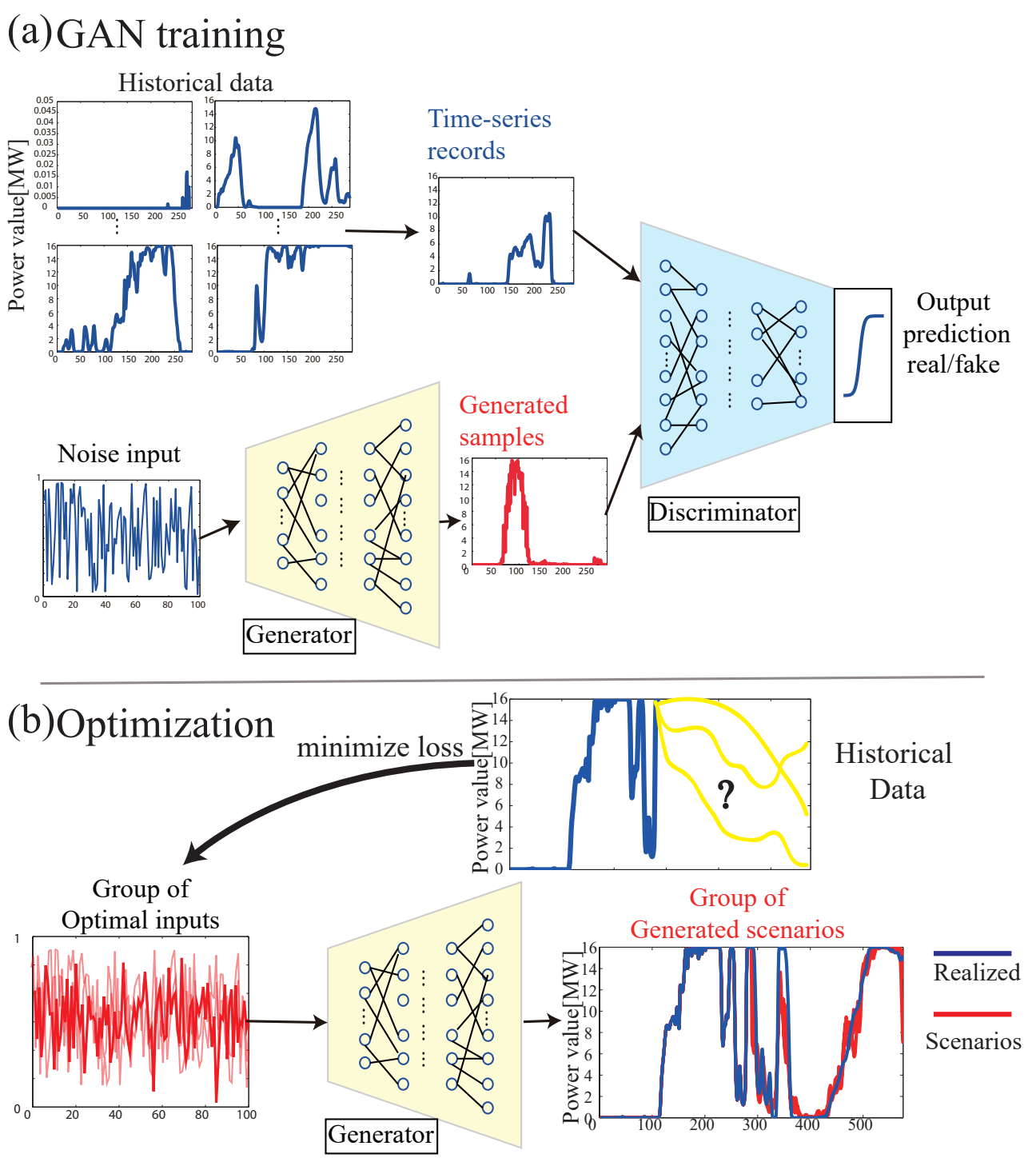

Figure 1. Model architecture for our proposed method. During the training process (a), the generator transforms noise vectors into generated time-series, and the discriminator is fed with both historical time-series and generated time-series, and try to discriminate the source of the input. Once training completes, we are able to optimize over the noise vectors to find the future scenarios from generator output conditioned on historical observations (b).

In [13], an unsupervised machine learning algorithm using Generative Adversarial Networks (GANs) [14] was introduced to directly generate realistic scenarios based only on historical data, without the need to fit an explicit model. In this paper, we extend the algorithm to the scenario forecasting problem. Compared to the work in [13], *this paper focuses on generating scenarios conditioned on a given central forecast*. Our proposed method is entirely *model-free and data driven*. Based on deep learning, the GANs used in our proposed method are unsupervised learners who can directly learn and generate time-series which hold same properties as the training time-series. The following optimization step would help us find group of scenarios based on point forecasts. Our approach can be used for a variety of scenario forecasts problems, e.g., wind and solar generation, and is easy to adjust the forecast horizon (e.g. ranging from 1-2 hours to 1-2 days) with little tuning. Specifically, we make the following contributions:

1) Based on any provided point forecast method along with historical observations, our method is able to generate a group of short-term forecasting scenarios representing the temporal correlations and fluctuation distribution.
2) The proposed approach can generate scenarios without forecast horizon or number restrictions.
3) The proposed approach is free of statistical assumptions and ready to use in real power generation processes.

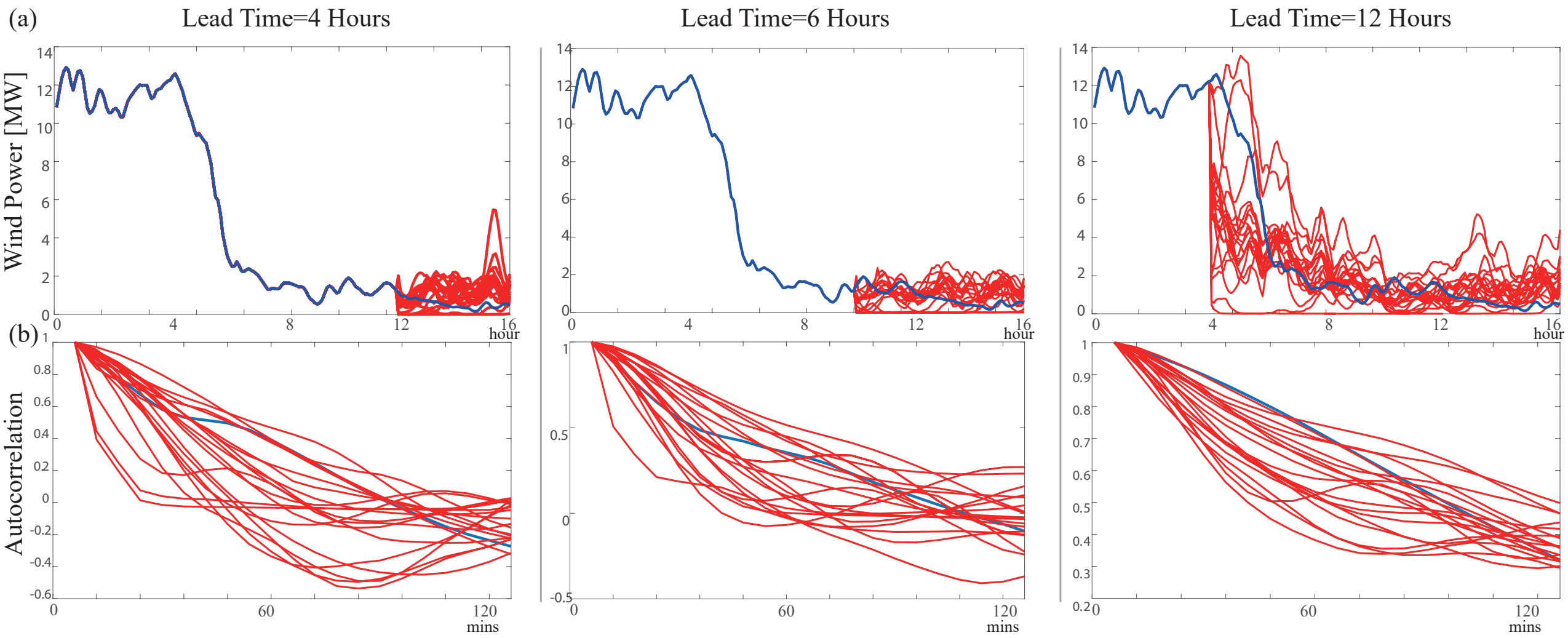

Figure 2. In (a) we show a group of 20 generated scenarios (red) with different scale of prediction lead time corresponding to the true measurements (blue). In (b) we show the associated autocorrelation plots for both measured values and generated scenarios.

Our proposed method for scenarios forecasts contains the following two components:

*1) Generative Adversarial Networks:* A generative adversarial networks (GANs) is composed of two deep neural networks, the generator and the discriminator, who play a zero-sum game. GANs are provided with past observations of renewables generation processes. Suppose these samples are drawn from an unknown underlying "true" distribution. The generator has access to a well-defined noise distribution (e.g., Gaussian) and can draw i.i.d. samples from this distribution. The generator's goal is to find a function that transforms a vector from the known noise distribution to a sample following the same distribution as past observations. The discriminator's goal is to distinguish the generated data and the true historical data. By training the generator and the discriminator to an equilibrium, the discriminator can no longer distinguish between generated and historical data, which means the generator can produce realistic time-series samples as if they are coming from the true distribution [14], [15]. Section. II describes this approach in more detail.

*2) Optimization of scenarios forecasts:* We are interested in forecasting a group of scenarios which could inform system operators the possible future realizations of power generation process. In Section. III we detail the setup of the optimization problem with a pre-trained GANs model. We also show how the optimization problem can be solved iteratively to obtain high-quality scenario forecasts. Some generated scenarios with different forecast horizons are shown in Fig. 2.

Numerical simulation are performed and evaluated in Section.IV. We will show our generated scenarios not only satisfy the needs of reliability and accurateness as a forecasting tool, they also capture the temporal dynamics of power generation. We also make our code for the proposed method freely available[1], which can meet the needs for an efficient computation tool for generating reliable and accurate scenarios.

## II. Generative Adversarial Networks

In this section, we describe the setup for GANs [14]. We first formulate the training objectives for the discriminator and the generator respectively, and show GANs is a good fit to generate a potentially unlimited number of renewable power production time-series. In Section III we illustrate how this time-series producer can be served in an optimization problem to find desired scenarios forecasts.

The architecture of GANs we use is shown in Fig. 1a. Assume observations $x_j^t$ for times $t \in T$ of renewable power production are available for each power plant $j$, $j = 1, \ldots, N$. Denote the true distribution of the observation as $\mathbb{P}_X$, which is unknown and maybe difficult to model because of complex spatial and temporal correlations. Suppose we have access to a group of noise vector input $z$ under a known distribution

[1]https://github.com/chennnnnyize/Scenario-Forecasts-GAN

$Z \sim \mathbb{P}_Z$ that is easily sampled from (e.g., jointly Gaussian or uniform). Given a sample $z$ drawn from $\mathbb{P}_Z$, our objective is to find a function $G$ such that after transformation, $G(z)$ follows $\mathbb{P}_X$. This is accomplished by simultaneously training two deep neural networks: the generator network $G(z;\theta^{(G)})$ and the discriminator network $D(x;\theta^{(D)})$. Here, $\theta^{(G)}$ and $\theta^{(D)}$ denote the weights of two neural networks, respectively. For convenience, we sometimes suppress the symbol $\theta$.

**Generator:** During the training process, the generator is trained to take a batch of inputs from the noisy distribution $\mathbb{P}_Z$, and by taking a series of up-sampling operations by neurons of different functions, and to output realistic time-series samples. Ideally, they should appear as if drawn from $\mathbb{P}_X$. Therefore, after training finishes, the mapping $G(z;\theta^{(G)})$ should follow the true data distribution $\mathbb{P}_X$.

**Discriminator**: The discriminator is trained simultaneously with the generator. It takes input samples either coming from real historical data or coming from the generator. By taking a series of operations of down-sampling using another deep neural network, it outputs a continuous value $p_{real}$ that measures to what extent the input samples belong to $\mathbb{P}_X$. The discriminator can be expressed as $D(x;\theta^{(D)})$, where $x$ may come from $\mathbb{P}_X$ or $\mathbb{P}_Z$. The discriminator is trained to learn to distinguish between $\mathbb{P}_X$ from $\mathbb{P}_Z$, and thus to maximize the difference between $D(X)$ ($X$ from real data) and $D(G(Z))$.

With the objectives for the discriminator and the generator defined, we can now formulate loss function $L_G$ for the generator and $L_D$ for the discriminator to train to optimize the performance of them (i.e., update neural networks' weights based on the losses). A small $L_G$ reflects that $G(z)$ is as realistic as possible from the discriminator's perspective, e.g., the generated scenarios are "looking like" historical scenarios to the discriminator. Similarly, a small $L_D$ indicates discriminator is good at telling the difference between generated scenarios and historical scenarios, which means there is a large difference between $\mathbb{P}_G$ and $\mathbb{P}_X$. Following this guideline and the loss defined in [15], we define $L_D$ and $L_G$ as:

$$L_G = -\mathbb{E}_Z[D(G(Z))] \quad (1a)$$
$$L_D = -\mathbb{E}_X[D(X)] + \mathbb{E}_Z[D(G(Z))]. \quad (1b)$$

In the above, the expectations are taken as empirical averages based either on the historical data or on the generated data. Note the functions $D$ and $G$ are parametrized by the weights of two distinct deep neural networks.

We can now combine (1a) and (1b) to construct the minimax game value function $V(G,D)$ for these two players:

$$\min_G \max_D V(G,D) = \mathbb{E}_X[D(X)] - \mathbb{E}_Z[D(G(Z))] \quad (2)$$

where $V(G,D)$ is the negative of $L_D$.

During first few training iterations, $G$ just generates time-series samples $G(z)$ totally different from samples in $\mathbb{P}_X$, and after learning from those samples coming from $\mathbb{P}_X$, the discriminator is able to reject $G(z)$ with high confidence. In that case, $L_D$ is small, and $L_G$, $V(G,D)$ are both large. The generator gradually learns to generate more realistically looking samples, while at the same time the discriminator is also trained to distinguish these newly fed generated samples from $G$. As training moves on and moves close the the equilibrium, $G$ is able to generate samples that look as realistic as real power generation time-series corresponding to a small $L_G$ value, while $D$ is unable to distinguish $G(z)$ from $\mathbb{P}_X$ with large $L_D$. Eventually, we are able to learn an unsupervised representation of the probability distribution of renewables time-series. By sampling $z$ from known distribution, we get $G(z)$ that appears "as if" it was sampled from the true distribution.

More formally, the minimax objective (2) of the game can be interpreted as the dual of the so-called Wasserstein distance (Earth-Mover distance) [16]. The Wasserstein distance between two distribution $\mathbb{X}$ and $\mathbb{Y}$ measures the effort (or "cost") needed to transport $\mathbb{X}$ to $\mathbb{Y}$. It is shown in [15] that we are precisely trying to get two distributions, $\mathbb{P}_x$ and $\mathbb{P}_{G(z)}$ to be close to each other by defined loss for $G$ and $D$ in (1a) and (1b) respectively.

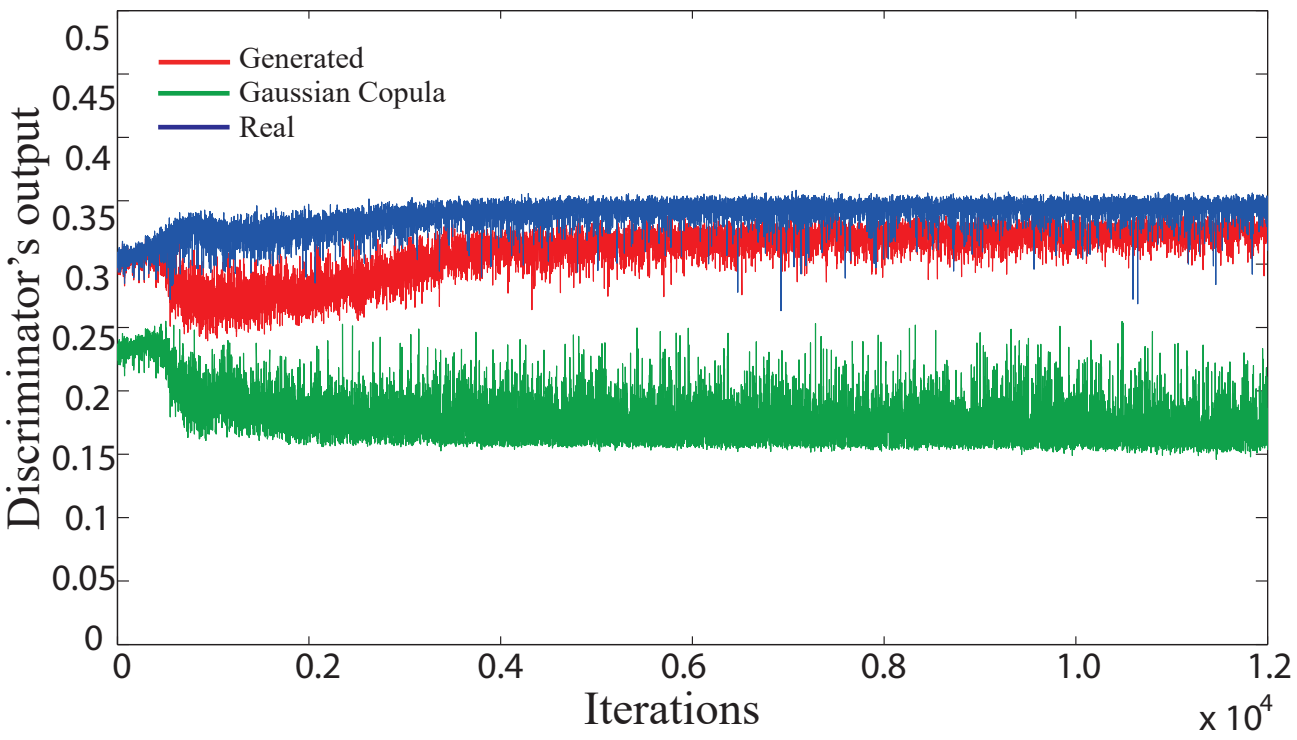

Figure 3. The evolution of discriminator $D$'s output during training process. Generated time-series(red) are indistinguishable from real measurements (blue). As a comparison, we also feed a batch of scenarios(green) generated by Gaussian Copula method, which is easy to distinguish from the discriminator's perspective.

Note that unlike previous approaches for generating scenarios given historical observations, which all involve the modeling of renewables generation stochastic processes [2], [10], [11], by using GANs we bypass the step of learning or modeling $\mathbb{P}_x$ explicitly. It also avoids complicated latent modeling in Variational Autoencoders [17]. The training algorithm of GANs for generating renewables time-series is summarized in Algorithm 1. In Fig. 3, we show the evolution of the loss function of the discriminator during training process and how the generated samples learn to mimic real historical data. In comparison, we show that the discriminator can always reliably detect that the scenarios generated using a Gaussian copula method [2] from the true realizations.

## III. Scenarios Forecasts Using GANs

In this Section, we show by formulating the scenario forecasts as an optimization problem, a trained GANs can be used

**Algorithm 1** GANs for Time-Series Generation

**Input:** Learning rate $\eta$, clipping parameter $c$, batch size $m$, Number of iterations for discriminator per generator iteration $n_{discri}$

**Initialize:** Initial weights $\theta^{(D)}$, $\theta^{(G)}$

**while** $\theta^{(D)}$ has not converged **do**

  **for** $t = 0, ..., n_{discri}$ **do**

    *# Update parameter for Discriminator*

    Sample batch from historical data:

    $\{x^{(i)}\}_{i=1}^{m} \sim \mathbb{P}_X$

    Sample batch from Uniform distribution:

    $\{z^{(i)}\}_{i=1}^{m} \sim Unif(-1,1)$

    Update discriminator nets using gradient descent:

    $g_{\theta^{(D)}} \leftarrow \nabla_{\theta^{(D)}}[-\frac{1}{m}\sum_{i=1}^{m} D(x^{(i)}) + \frac{1}{m}\sum_{i=1}^{m} D(G(z^{(i)}))]$

    $\theta^{(D)} \leftarrow \theta^{(D)} - \eta \cdot RMSProp(\theta^{(D)}, g_{\theta^{(D)}})$

    $\theta^{(D)} \leftarrow clip(\theta^{(D)}, -c, c)$

  **end for**

  *# Update parameter for Generator*

  Update generator nets using gradient descent:

  $g_{\theta^{(G)}} \leftarrow \nabla_{\theta^{(G)}} \frac{1}{m}\sum_{i=1}^{m} D(G(z^{(i)}))$

  $\theta^{(G)} \leftarrow \theta^{(G)} - \eta \cdot RMSProp(\theta^{(G)}, g_{\theta^{(G)}})$

**end while**

to generate a group of scenarios given past observations and point forecasts.

### A. Mathematical Formulation

For a typical renewable power generation site, assume at timestep $t$, we have records for actual past power outputs $\mathbf{p}_{hist} \in \mathbf{R}^{h+1}$ with $\mathbf{p}_{hist} = [p_t \; \dots \; p_{t-h}]^T$. Meanwhile, we have some forecasting method to obtain the point forecasts $\hat{p}_{t+i}$ for each look-ahead time $i = 1, ..., k$ given $\mathbf{p}_{hist}$. This forecast is denoted by $\hat{\mathbf{p}}_{pred} = [p_{t+1} \; \dots \; p_{t+k}]^T$, where $k$ is the forecasting horizon. Based on the historical information and the point forecast, we are interested in generating a group of $N$ scenarios $S = \{s_1, ..., s_N\}$, $s_i \in \mathbf{R}^k$, which represent the possible variations around the point forecast and accurately reflect the temporal dynamics of future generation. Note we focus on the scenario forecasting problem, so the central point forecast can be provided by any method.

Assume we have trained a GANs model based on the set of observations. Given some input noise $z$, $G(z)$ generates a possible realization without regarding to the historically observed data $\mathbf{p}_{hist}$ and the point forecast $\hat{\mathbf{p}}_{pred}$. Therefore, we need to *constrain the possible $G(z)$ to satisfy two conditions:* 1) the part of $G(z)$ from time index $t-h$ to $t$ should be close to the historical data $\mathbf{p}_{hist}$; 2) the part of $G(z)$ from time index $t+1$ to $t+k$ should be realistic and respect the point forecast.

To describe these conditions, we introduce two projection operators that separate a vector into two parts. Given $\mathbf{v} \in \mathbf{R}^{h+k+1}$, we denote two projection operations $\mathbb{P}_{hist}(f(\mathbf{v})) = [v_1, ..., v_{h+1}]^T$ and $\mathbb{P}_{pred}(f(\mathbf{v})) = [v_{h+2}, ..., v_{h+k+1}]^T$ to extract former $h+1$ and latter $k$ dimensions of $\mathbf{v}$, respectively.

Meanwhile, we want to constrain generated scenarios do not conflict with the information provided by point forecasts $\hat{\mathbf{p}}_{pred}$ (e.g., information from numerical weather prediction (NWP)). Then we can constrain latter part of $G(z)$ so that they do not conflict with the given $\hat{\mathbf{p}}_{pred}$. Then to ensure that the first part of $G(z)$ resembles $\mathbb{P}_{hist}$, we use the following cost function:

$$||\mathbb{P}_{hist}(G(z)) - \mathbf{p}_{hist}||_2. \tag{3}$$

To ensure the generated scenarios are realistic, we add a loss term $-D(G(z))$ where $D$ is the discriminator output (recall larger discriminator ouput indicates more realistic samples). Finally, we use the point forecast $\hat{\mathbf{p}}_{pred}$ by defining a *prediction interval* that the generated scenarios should lie in [18]. We describe this interval with an upper bound $U^{\alpha}(\hat{\mathbf{p}}_{pred})$ and a lower bound $L^{\alpha}(\hat{\mathbf{p}}_{pred})$, controlled by a parameter $\alpha$ (can be interpreted as the prediction confidence):

$$L^{\alpha}(\hat{\mathbf{p}}_{pred}) = \frac{1}{\alpha}\hat{\mathbf{p}}_{pred}, \quad U^{\alpha}(\hat{\mathbf{p}}_{pred}) = \alpha\,\hat{\mathbf{p}}_{pred} \tag{4}$$

Using all of the objectives and constraints above, given the observation and forecast vector pair $\mathbf{p}_{hist}$, $\hat{\mathbf{p}}_{pred}$, and GANs pre-trained model $G$, $D$, the scenario forecasts problem can be formulated as a constrained optimization problem:

$$\min_{z} \quad ||\mathbb{P}_{hist}(G(z)) - \mathbf{p}_{hist}||_2 - \gamma\, D(G(z)) \tag{5}$$

$$s.t. \quad z \in \mathcal{Z} \tag{5a}$$

$$L^{\alpha}(\hat{\mathbf{p}}_{pred}) \leq \mathbb{P}_{pred}(G(z)) \leq U^{\alpha}(\hat{\mathbf{p}}_{pred}). \tag{5b}$$

where $\gamma$ is a weighting parameter; (5a) constrains $z$ to be within the domain of $G(z)$, which we take to be a hypercube $\mathcal{Z} = [-1,1]^{h+k+1}$; (5b) constrains the generators' output to be within the given prediction intervals given $\hat{\mathbf{p}}_{pred}$. By solving above optimization problem, we can obtain a forecasting scenario $\mathbb{P}_{pred}(G(z^*))$.

Since both of the objective and constraints in (5) are non-convex, to deal with the inequality constraints (5b), we propose to substitute it into the main objective with two log barriers. Then the optimization problem is reformulated as

$$\min_{z} \quad ||\mathbb{P}_{hist}(G(z)) - \mathbf{p}_{hist}||_2 - \beta\; \log(\mathbb{P}_{pred}(G(z)) - L^{\alpha}(\hat{\mathbf{p}}_{pred})) - \beta\; \log(U^{\alpha}(\hat{\mathbf{p}}_{pred}) - \mathbb{P}_{pred}(G(z))) - \gamma\, D(G(z)) \tag{6a}$$

$$s.t. \quad z \in \mathcal{Z}. \tag{6b}$$

where $\beta$ is the weighting parameter for log barriers.

In next subsection, we will illustrate how we are able to find a *group of solutions $z$* for problem defined in (6) given the fixed, differentiable, yet non-convex function $G(\cdot)$.

### B. Forecasting Scenarios with Pre-trained GANs

Because of the highly non-convex nature of $G(\cdot)$ and $D(\cdot)$, there exist many local optima in (5). The key to finding a group of solutions to (5) exploits that fact. Figure 4 shows the landscape of solutions in a one-dimensional illustration. To ensure that we reach a good local optimum, we add momentum to the gradient descents algorithm [19] to skip from saddle points and shallow local optima.[2]

[2]There is a growing body of literature on the local optima of non-convex functions and interested readers can refer to [19] and the references within.

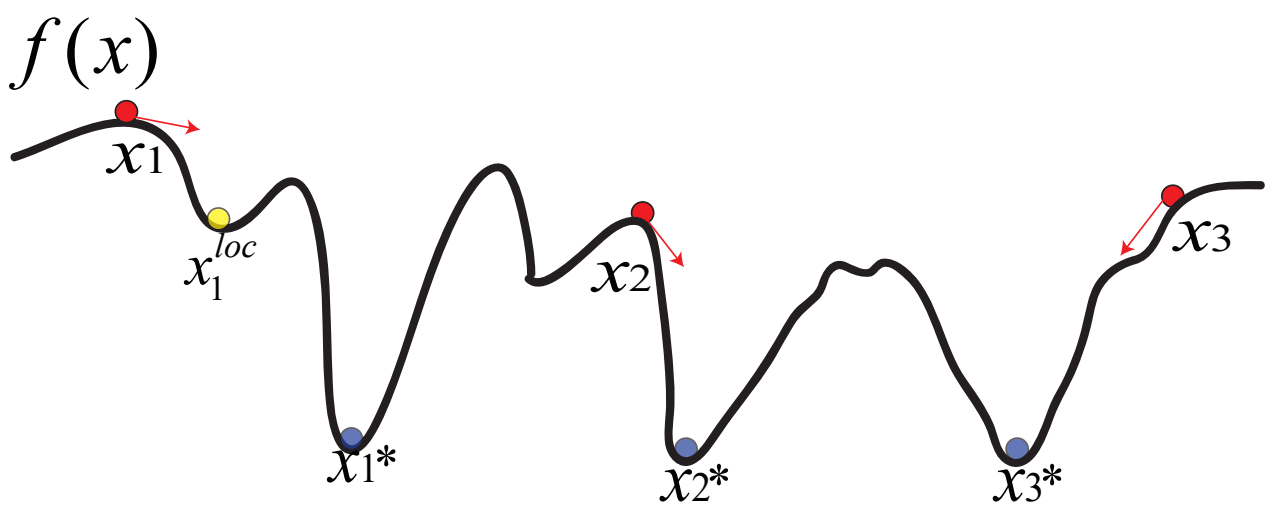

Figure 4. Several optimal points exist on the 1-D curve, while different initial points and momentum gradient descents can help reach these optimal points by getting over the local optimal points.

Since there are multiple local optima to (5), we can start at different initial points $z_i \in \mathcal{Z}$ and find distinct forecasting scenario $\mathbb{P}_{pred}(G(z_i^*))$ by solving (6) using gradient descents with momentum (MomentumGD). As the training loss defined in (2) incurs $G$ to generate diverse modes given different $z$, we are able to obtain a group of distinct yet realistic scenarios with different initial starting values.

In order to obtain good starting points for $z$ which $\mathbb{P}_{pred}(G(z)))$ do not fall outside of the log barriers in (6), we first solve the following subsidiary problem:

$$\min_{z} \quad ||\mathbb{P}_{pred}(G(z)) - \mathbf{p}_{initial}||_2 \tag{7}$$

$$s.t. \quad z \in \mathcal{Z} \tag{7a}$$

$$L^{\alpha}(\hat{\mathbf{p}}_{pred}) \leq \mathbf{p}_{initial} \leq U^{\alpha}(\hat{\mathbf{p}}_{pred}). \tag{7b}$$

where $\mathbf{p}_{initial}$ is sampled uniformly at random from $[L^{\alpha}(\hat{\mathbf{p}}_{pred}),\ U^{\alpha}(\hat{\mathbf{p}}_{pred})]$. In order that we can always obtain a good initial $z$, we set $\alpha$ in (7) to be slightly smaller than $\alpha$ in main objective function (6). In Algorithm 2 we summarize our approach for generating a group of scenarios provided with a pre-trained GANs weights as well as pairing historical measurements and point forecasts.

## IV. Simulation Results

In this section, we study the performance of the proposed method for scenario forecasts over various forecasting horizons. We focus on wind power production, and show the scenarios generated by our method not only conform to the statistical properties of real measurements, but also capture the spatial and temporal correlations. We also show that our approach is flexible to generate scenarios for varying prediction intervals as well as prediction horizons. All experiments are implemented using Python 2.7 with deep learning open-source package TensorFlow [20]. The GANs training procedure is accelerated by two Nvidia Geforce GTX TITAN X GPUs.

Both the generator and the discriminator deep neural networks are composed of two convolutional layers and two fully-connected layers. All models in this paper are trained using RmsProp optimizer [21], which is a learning-rate self-adaptive gradient descent algorithm. Weights for neurons in both neural networks were initialized from a centered normal distribution with standard deviation of 0.02. Batch normalization is adopted before every layer except the input layer to stabilize learning by normalizing the input of every layer

**Algorithm 2** Forecasting Scenarios with GANs

**Input:** Pre-trained GANs model weights $\theta^{(G)}$, $\theta^{(D)}$
**Input:** Measurements $\mathbf{p}_{hist}$, point forecast $\hat{\mathbf{p}}_{pred}$
**Input:** PI level $\alpha$, initial PI level $\alpha_{sub}$, weighting parameters $\gamma, \beta$, initial point finder iterations $n_{init}$, scenario finder iterations $n_{scen}$, learning rate $\eta$, scenario number $N$
**Initialize:** Generated scenarios S $\leftarrow \emptyset$
  **for** $scenario = 0, \ldots, N$ **do**
    Sample $\mathbf{p}_{initial} \sim Unif(L^{\alpha_{sub}}(\hat{\mathbf{p}}_{pred}), U^{\alpha_{sub}}(\hat{\mathbf{p}}_{pred}))$
    Sample $z \sim Unif(-1, 1)$
    *# Find good initial z*
    **for** $iteration = 0, \ldots, n_{init}$ **do**
      Update $z$ using gradient descent:
      $g_z \leftarrow \nabla_z L_{sub}$ *#*$L_{sub}$ *is defined by 7*
      $z \leftarrow z - z \cdot MomentumGD(z, g_z)$
      $z \leftarrow clip(z, -1, 1)$
    **end for**
    *# Find forecasting scenarios*
    **for** $iteration = 0, \ldots, n_{scen}$ **do**
      Update $z$ using gradient descent:
      $g_z \leftarrow \nabla_z L_{main}$ *#*$L_{main}$ *is defined by 6*
      $z \leftarrow \theta^{(D)} - \eta \cdot MomentumGD(z, g_z)$
      $z \leftarrow clip(z, -1, 1)$
    **end for**
    S.insert(G(z))
  **end for**

to have zero mean and unit variance. With exception of the output layer, rectified linear units (ReLU) has been used as the activation function in the generator, while Leaky-ReLU is used in the discriminator. We observed in Algorithm 1, a setting of $n_{discri} = 4$ could get the fastest convergence rate for $D$. Once the discriminator has converged to similar outputs value for $D(G(Z))$ and $D(X)$, the generator was able to generate realistic power generation samples.

### A. Description of Data

In order to test the performance of our proposed framework for scenario generation, we set up our numerical simulations based on wind power data published by the NREL Wind Integration National Data Set (WIND) Toolkit [3]. Actual power measurements have a resolution of 5 minutes. The dataset also contains deterministic, day-ahead forecasts along with estimated 10% and 90% forecast quantiles. The detailed NWP-based forecasts method is described in [22]. We construct our dataset by aggregating 20 wind turbines' records from Jan.1st, 2007 to Dec. 31st, 2013. All these wind turbines are located in WA, USA and are of geographical proximity. Selected wind turbines have a nominal capacity of $16MW$. In total there are $14,728,320$ measurements, and we split 90% of daily samples as the real training data for our GANs model, while the remaining 10% samples are only used to test the performance of the proposed scenario forecasts method. All wind power measurements and forecasts are normalized to $[0, 1]$.

[3]https://www.nrel.gov/grid/wind-integration-data.html

### B. Validation Framework

The validation of the quality of generated scenarios is more complex than evaluating the performance of point forecasts. On the one hand, generated scenarios should be realistic enough to reflect the interdependence structure of forecasting values at different prediction horizons; on the other hand they should represent all the possible future realizations given past observations at certain wind farm.

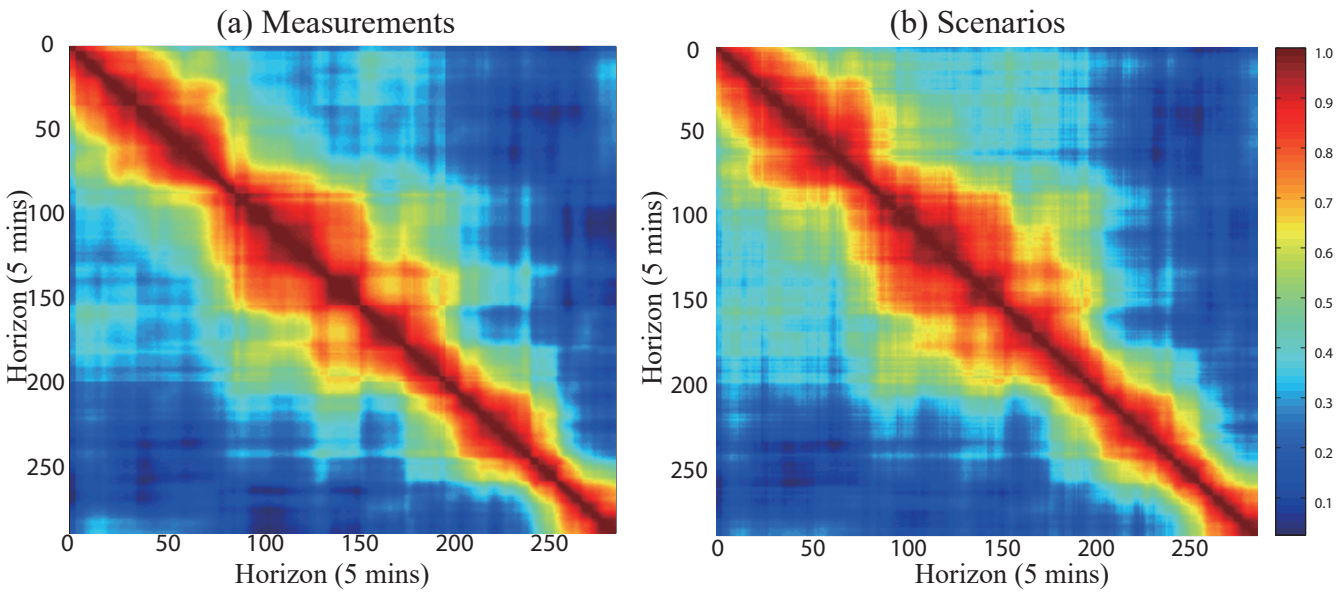

Figure 5. The Pearsons correlation matrix for a group of historical realizations (left) and scenarios forecasts (right) indicates our scenario forecasts method capture the linear correlation for forcasting lead time varying from 5 minutes to 1 day.

First, we examine generated scenarios' temporal statistics. We calculate and compare samples' autocorrelation with respect to look-ahead time $k$:

$$R(k) = \frac{E[(s_t - \mu)(s_{t+k} - \mu)]}{\sigma^2} \quad (8)$$

where $s \in S$ represents sample either of generated scenarios or realizations with mean $\mu$ and variance $\sigma$.

We also make use of the Pearson's correlation coefficient, which is a standard method to evaluate the linear relationship of time-series at various look-ahead times. Given the set of generated scenarios or realizations $S$, each term $\rho_{i,j}$ in the Pearson's correlation matrix denotes the Pearson correlation for lead time $i$ and $j$, and is calculated by

$$\rho_{i,j} = \frac{Cov(S_i, S_j)}{\sigma_{S_i}\sigma_{S_j}} \quad (9)$$

where $Cov(S_i, S_j)$ is the covariance of $S_i$ and $S_j$.

In order to verify the group of generated scenarios are able to represent possible future realizations, the scenarios should be able to cover the actual value of power generation (reliable), while at the same time distance between generated scenarios should be small (sharp). We make use of the Continuous Ranked Probability Score (CPRS) [23], which is a negatively-oriented score (smaller scores are better). It is a comprehensive metric that jointly evaluates the reliability and sharpness of generated scenarios. The score at lead time $k$ is defined as

$$CPRS_k = \frac{1}{N}\sum_{t=1}^{N}\int_0^1 (\hat{F}_{t+k|t}(p) - I(p \geq p_{t+k}))^2 dp \quad (10)$$

where $N$ is the total number of evaluated scenarios, $\hat{F}_{t+k|t}(p)$ is the cumulative distribution for normalized generated scenarios' value at lead time $k$, and $I(p \geq p_{t+k})$ is the indicator function to compare the normalized scenarios and measurements. Since we are not using quantile statistics to calculate $\hat{F}_{t+k|t}(p)$, we use the discrete-valued $\hat{F}_{t+k|t}(p)$ to calculate (10).

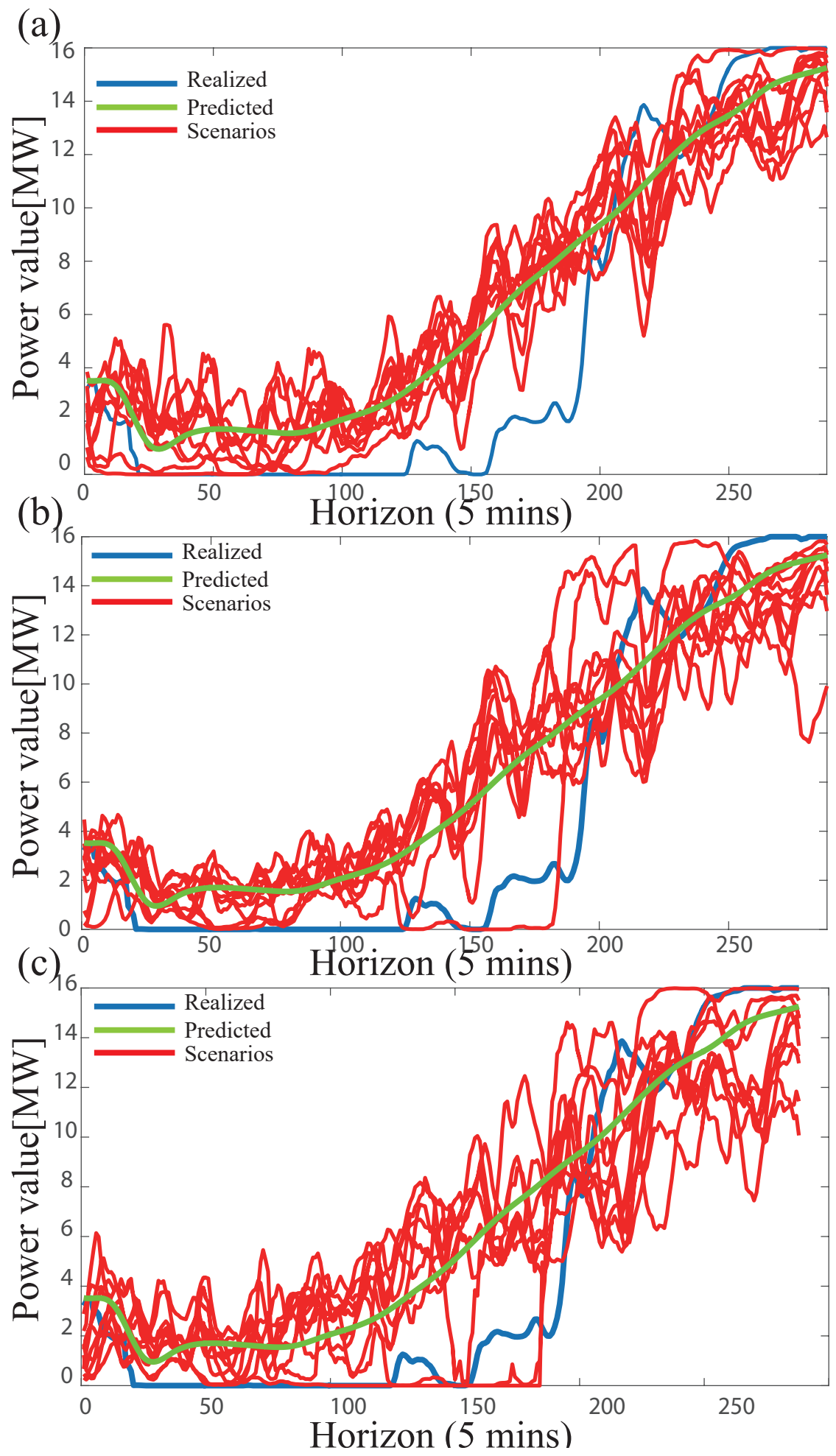

Figure 6. Plot (a) (b) (c) correspond to group of 10 day-ahead scenarios (red) with varying PIs of 1.5, 2 and 3 respectively.

### C. Simulation Results

Recall that Fig. 3 showed the output evolution of $D$, where $D(x)$ and $D(G(z))$ converged after about 6,000 training iterations. In addition, we evaluate the quality of a generated time-series from empirical Gaussian Copula method [2]. In this case, $D$ is able to distinguish the generated samples from real measurements. This observation suggests that eventhough Gaussian Copula method tries to model the interdependence structure for time-series, the generated scenarios are still different from actual realizations.

We then validate if scenarios coming from our method have similar temporal correlation as the actual wind power values. In Fig. 5 we plot the colormap for the covariance matrix of a group of 32 wind turbines' 24−hour actual measurements, along with 1,600 forecasting scenarios with 50 scenarios for each realization. The $x-$ and $y-$ axes are for the prediction horizon $k$. Similar covariance matrix element values indicate that without any model assumptions being made during training process, our proposed scenario generation method is able to capture the temporal dependency accurately.

Fig. 2 shows a group of 20 generated scenarios with forecasts lead time ranging from $4h$ to $12h$. We show that by only changing projection length $k$, our approach is able to conveniently generate reliable and sharp scenarios for different forecast horizons. Meanwhile, these scenarios' autocorrelation plots cover the realizations, which indicate the generated scenarios are able to represent the temporal dependence of any length.

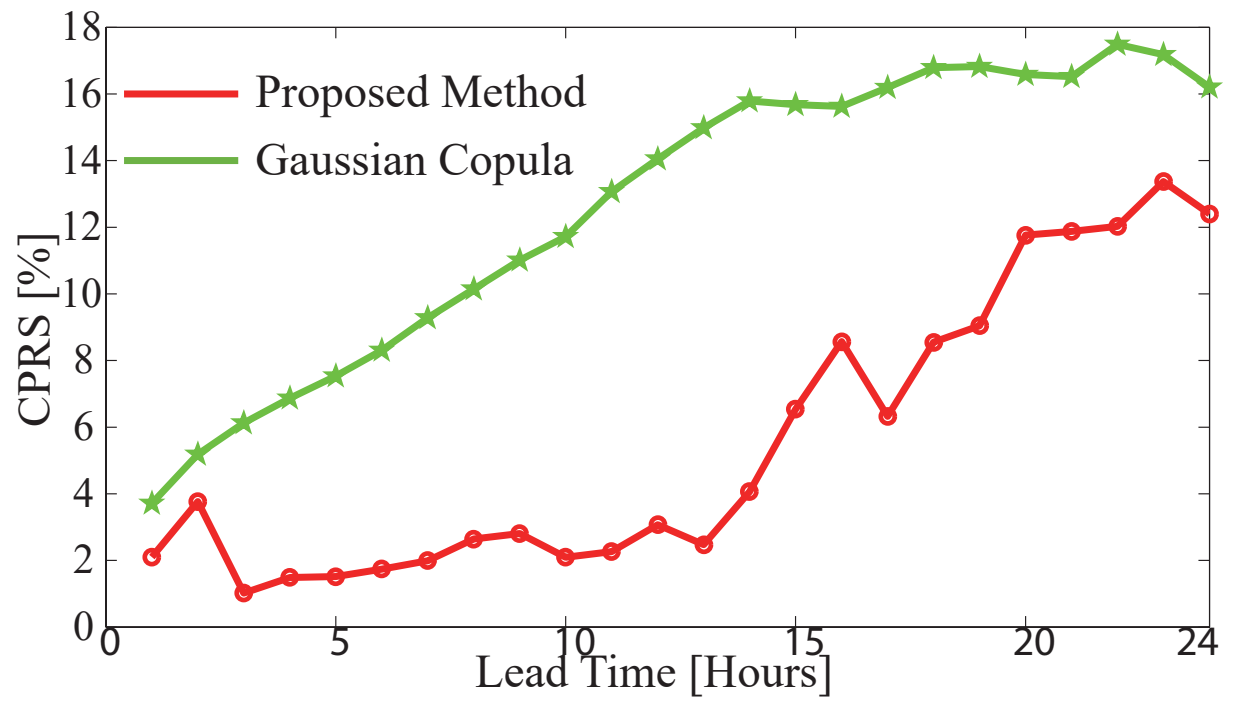

Figure 7. CPRS of scenarios generated by proposed method and empirical Gaussian Copula [2].

In Fig. 6 we specifically select one 24−hour sample whose point forecast is deviating a lot from the actual measurements. By selecting different prediction intervals, our proposed method could reflect the trade-off between reliability and sharpness. When the interval level is $\alpha = 1.5$, generated scenarios are close to point forecasts, yet fail to cover the realizations; while when $\alpha = 3$, generated scenarios could cover the actual power production values, but are less concentrated.

The performance of the proposed method is also demonstrated by the CPRS score. Results for our approach and Gaussian copula method are plotted in Fig. 7. Both approaches use the same training dataset to get the time-series generator or to find the estimate of the covariance matrix, and are tested on the stand-alone testing samples. The proposed method has better performance at different lead time compared to Gaussian Copula. Since point forecasts normally accumulate larger errors with longer forecast horizons, both methods have growing CPRS values with respect to forecasting horizons.

## V. Conclusion

In this paper we proposed a data-driven unsupervised machine learning approach for forecasting scenarios of renewables power generation processes. The proposed method is flexible and easily implemented in problems with high penetration of renewables. Numerical results show that comparing with existing scenario generation approaches, the proposed method is able to generate realistic, high quality scenarios capturing spatiotemporal behaviors of renewables without any explicit model construction.